\documentclass[11pt]{article}

\usepackage{amssymb,latexsym,amsmath}
\begin{document}

\newcommand{\bfi}{\bfseries\itshape}

\makeatletter

\@addtoreset{figure}{section}

\def\thefigure{\thesection.\@arabic\c@figure}

\def\fps@figure{h, t}

\@addtoreset{table}{bsection}

\def\thetable{\thesection.\@arabic\c@table}

\def\fps@table{h, t}

\@addtoreset{equation}{section}

\def\theequation{\thesubsection.\arabic{equation}}

\makeatother

\newcommand{\ea}{\mbox{{\bf a}}}
\newcommand{\eu}{\mbox{{\bf u}}}
\newcommand{\ueu}{\underline{\eu}}
\newcommand{\oeu}{\overline{\eu}}
\newcommand{\ew}{\mbox{{\bf w}}}
\newcommand{\ef}{\mbox{{\bf f}}}
\newcommand{\eF}{\mbox{{\bf F}}}
\newcommand{\eC}{\mbox{{\bf C}}}
\newcommand{\en}{\mbox{{\bf n}}}
\newcommand{\eT}{\mbox{{\bf T}}}
\newcommand{\eL}{\mbox{{\bf L}}}
\newcommand{\eV}{\mbox{{\bf V}}}
\newcommand{\eU}{\mbox{{\bf U}}}
\newcommand{\ev}{\mbox{{\bf v}}}
\newcommand{\uev}{\underline{\ev}}
\newcommand{\eY}{\mbox{{\bf Y}}}
\newcommand{\eK}{\mbox{{\bf K}}}
\newcommand{\eP}{\mbox{{\bf P}}}
\newcommand{\eS}{\mbox{{\bf S}}}
\newcommand{\eJ}{\mbox{{\bf J}}}
\newcommand{\eB}{\mbox{{\bf B}}}
\newcommand{\leb}{\mathcal{ L}^{n}}
\newcommand{\eI}{\mathcal{ I}}
\newcommand{\eE}{\mathcal{ E}}
\newcommand{\hen}{\mathcal{ H}^{n-1}}
\newcommand{\eBV}{\mbox{{\bf BV}}}
\newcommand{\eSBV}{\mbox{{\bf SBV}}}
\newcommand{\eBD}{\mbox{{\bf BD}}}
\newcommand{\eSBD}{\mbox{{\bf SBD}}}
\newcommand{\ecs}{\mbox{{\bf X}}}
\newcommand{\eg}{\mbox{{\bf g}}}
\newcommand{\paromega}{\partial \Omega}
\newcommand{\gau}{\Gamma_{u}}
\newcommand{\gaf}{\Gamma_{f}}
\newcommand{\sig}{{\bf \sigma}}
\newcommand{\gac}{\Gamma_{\mbox{{\bf c}}}}
\newcommand{\deu}{\dot{\eu}}
\newcommand{\dueu}{\underline{\deu}}
\newcommand{\dev}{\dot{\ev}}
\newcommand{\duev}{\underline{\dev}}
\newcommand{\weak}{\rightharpoonup}
\newcommand{\weakdown}{\rightharpoondown}

\title{PREPRINT SERIES OF THE INSTITUTE OF MATHEMATICS \\  
OF THE ROMANIAN ACADEMY \\  $\left. \right.$\\  Perturbed area functionals and brittle damage mechanics}

\author{Marius Buliga \\
 \\
Institute of Mathematics, Romanian Academy \\
P.O. BOX 1-764, RO 70700\\
Bucure\c sti, Romania\\
{\footnotesize Marius.Buliga@imar.ro}}

\date{Preprint No. 27/1996}

\maketitle

\begin{abstract}
Some Mumford-Shah functionals are revisited as perturbed area functionals in connection with brittle 
damage mechanics. We find minimizers "on paper" for the classical Mumford-Shah functional for some particular two dimensional domains and boundary conditions. These solutions raise the possibility of validating experimentally the energetic model of crack appearance. Two models of brittle damage and fracture are proposed after; in the one of these models the crack belongs to the set of integral varifolds. We have felt the necessity to start the paper with a preliminary section  concerning classical results in equilibrium of a cracked 
elastic body reviewed in the context of Sobolev spaces with respect to a measure. 
\end{abstract}

\thispagestyle{empty}

\newpage

\section{Preliminaries}

\subsection{Functions with bounded variation or deformation}
\indent

Let $\Omega \subset R^{n}$ be an open set  with local Lipschitzian boundary. $\leb$ or $\mid \cdot \mid$ denotes the Lebesgue measure on $R^{n}$ ; $\hen$ is the 
Hausdorff $n-1$-dimensional measure and $\eB(\Omega)$ denotes the collection of Borelian sets from $\Omega$ or the family of Borelian maps from $\Omega$ to $R$.
For any $\eu \in L^{1}(\Omega , R^{m})$ let $$\tilde{\eu}: \Omega_{\eu} \rightarrow R^{m}$$ denote 
the exact representative of $\eu$ ; at each point $x$ from the Lebesgue set $\Omega_{\eu}$ of $\eu$ $\tilde{\eu}(x)$ is the approximate limit of $\eu$ in $x$. The complementary 
set of $\Omega_{\eu}$ is $\eS_{\eu}$. $\eS_{\eu}$ is a Borelian set with Lebesgue measure zero. For a given vector $\nu \in \eS^{n-1}$ one can define the approximate limit 
in the $\nu$ direction in a point $x$ as the number $\tilde{\eu}^{\nu}(x)$ which satisfies : $$\lim_{\rho \rightarrow 0} \frac{1}{\rho^{n}} \int_{ \left\{ y \in B_{\rho}(x) \mbox{ : } 
\langle y-x , \nu \rangle > 0 \right\} } \mid \eu (y) - \tilde{\eu}^{\nu}(x) \mid \mbox{ d}x = 0$$ If $x \in \Omega_{\eu}$ then  $\tilde{\eu}^{\nu}(x)$ exists for any $\nu$ and is equal 
to $\tilde{\eu}(x)$; denote then by $\eJ_{\eu} \subset \eS_{\eu}$ the set of jump points of $\eu$ , i.e. the set of points  $x \in \eS_{\eu}$ for which $\tilde{\eu}^{\nu}(x)$ and 
$\tilde{\eu}^{-\nu}(x)$ both exists for a point-dependent $\nu$.

The space $\eBV(\Omega,R^{m})$ of $R^{m}$-valued functions of bounded variation on the open set $\Omega$ is the subset of $L^{1}(\Omega,R^{m})$ of all functions whose distributional 
derivative, regarded as a set measure, has bounded variation. In the same way one can define in the natural way the space $\eBV_{loc}(R^{n},R^{m})$ of functions with locally bounded variation  
which is a subset of $L^{1}_{loc}(R^{n},R^{m})$. $\eBV(\Omega,R^{m})$ is a Banach space endowed with the norm: $$ \| \eu \|_{BV} = \| \eu \|_{L^{1}} + \mid D\eu \mid (\Omega)$$ where
$\mid D\eu \mid (B)$ is the variation of $D\eu$ over $B \in \eB(\Omega)$. 

The space $\eBD(\Omega)$ of $R^{n}$-valued functions of bounded deformation on $\Omega$ is the subset of $L^{1}(\Omega,R^{m})$ of all functions whose symmetric distributional 
derivative $E\eu$ has bounded variation. It is a Banach space too endowed with the norm $$ \| \eu \|_{BD} = \| \eu \|_{L^{1}} + \mid E\eu \mid (\Omega)$$ Obviously 
$\eBV(\Omega,R^{n}) \subset \eBD(\Omega)$. The resemblances and the differences between these two spaces are visible in the following results.

\vspace{.5cm} 

{\bf Theorem 1.1.} {\it Let $\eu \in L^{1}(\Omega,R^{m})$. Then 
\begin{itemize}
\item (De Giorgi) If $\eu \in \eBV(\Omega,R^{m})$ then $\eS_{\eu}$ is countably rectifiable, $\hen(\eS_{\eu} \setminus \eJ_{\eu}) = 0$ and in $\hen$-almost every point  $x \in \eS_{\eu}$
exists the approximate limits of $\eu$ in the directions $\nu(x)$ and $-\nu(x)$ where $\nu(x)$ is the normal to $\eS_{\eu}$ in $x$. 
\item (Kohn, Ambrosio, Coscia, Dal Maso) Let $m=n$ and $\eu \in \eBD(\Omega)$. Let $\Theta_{\eu}$ be the Kohn set :
$$\Theta_{\eu} = \left\{ x \in \Omega \mbox{ : } \limsup_{\rho \rightarrow 0^{+}} \frac{ \mid E\eu \mid (B_{\rho}(x))}{\rho^{n-1}} \ > \ 0 \right\}$$
Then $\Theta_{\eu}$ is countably rectifiable , $\eJ_{\eu} \subseteq \Theta_{\eu}$ and $\hen(\Theta_{\eu} \setminus \eJ_{\eu}) = 0$~.
\end{itemize} }

\vspace{1.cm}

{\bf Theorem 1.2.}{ \it Let $\eu \in L^{1}(\Omega,R^{m})$. Then
\begin{itemize}
\item (Calderon, Zygmund) If \ $\eu \in \eBV(\Omega,R^{m})$ \ then \ $\eu$ \ is approximately \ differentiable \newline $\leb$-a.e. in $\Omega$. 
The approximate differential map $x \mapsto \nabla \eu (x)$ is integrable. 
$D\eu$ splits into three mutually singular measures on $\Omega$
$$D\eu = \nabla \eu \mbox{ d}x \ + \ [\eu] \otimes \nu \hen_{|_{\eS_{\eu}}} \ + \ C\eu$$
where $[\eu]$ is the jump of $\eu$ in respect with the normal direction on $\eS_{\eu}$ $\nu$. $C\eu$ is the Cantor part of $D\eu$ defined by $C\eu(A) = D^{s}\eu(A \setminus \eS_{\eu})$ where $D^{s}\eu$ is 
the singular part of $D\eu$ in respect to $\leb$.
\item (Ambrosio, Coscia, Dal Maso)  Let $m=n$ and $\eu \in \eBD(\Omega)$. Then $\eu$ has symmetric approximate differential $\epsilon(\eu)$ $\leb$-a.e. in $\Omega$ and $E\eu$ splits into three mutually
singular measures on $\Omega$
$$E\eu = \epsilon (\eu) \mbox{ d}x \ + \ [\eu] \odot \nu \hen_{|_{\eJ_{\eu}}} \ + \ E^{c}\eu$$ 
Moreover $\eu$ is approximately differentiable $\leb$-a.e. in $\Omega$.
\end{itemize} }

\vspace{1.cm}

The previous theorem allows the introduction of two new spaces. The space   of special functions with bounded variation $\eSBV(\Omega,R^{m})$ is the subset of $\eBV(\Omega,R^{m})$ of 
all functions whose Cantor part $C\eu$ is equal to zero. The space $\eSBD(\Omega)$ of special functions of bounded deformation 
is the subset of $\eBD(\Omega)$ of all functions whose Cantor part $E^{c}\eu$ is equal to zero. These special functions have the following regularity properties:

\vspace{.5cm}

{\bf Theorem 1.3} { \it  
\begin{itemize}
\item $W^{1,1}(\Omega,R^{m}) \subset \eBV(\Omega,R^{m})$. The inclusion is continuous in respect with the Banach space topologies. If $$\eu \in \eSBV(\Omega,R^{m})$$ then 
$$\eu \in W^{1,1}(\Omega \setminus \eS_{\eu},R^{m})$$ Moreover if $\eu \in W^{1,1}(\Omega \setminus K,R^{m})\cap L^{\infty}(\Omega , R^{m})$ , where $K$ is a closed , 
countably rectifiable set with $\hen(K) < + \infty$, then $\eu \in \eSBV(\Omega,R^{m})$ and $\hen(K \setminus \eS_{\eu}) = 0$.
\item Let $LE^{1}(\Omega)$ be the Banach space of $L^{1}(\Omega,R^{n})$ functions with $L^{1}$ symmetric differential.  If $\eu \in \eSBD(\Omega)$ then $\eu \in LE^{1}(\Omega\setminus \eJ_{\eu})$.
Let $K$ be a closed , countably rectifiable set with $\hen(K) < + \infty$. If $\eu \in LE^{1}(\Omega\setminus K) \cap  L^{\infty}(\Omega , R^{n})$ then $\eu \in \eSBD(\Omega)$ and 
 $\hen(K \setminus \eJ_{\eu}) = 0$.
\end{itemize}  }

Remind also that if $\eu \in \eSBD(\Omega)$ then the extension of $\eu$ by $0$ in $R^{n} \setminus \overline{\Omega}$ is a $\eSBD_{loc}(R^{n})$ map.

\subsection{Sobolev spaces with respect to a measure and fractured media}
\indent

Let $\mu$ be a finite positive measure on $\Omega$, where $\Omega$ is an open bounded subset of 
$R^{n}$, with Lipschitz boundary. For any $p > 1$ 

$$W^{1,p}_{\mu} (\Omega, R^{k}) = \left\{ \eu \in \eBV(\Omega, R^{k}) \mbox{ : } \mid D\eu 
\mid \ll \mu \mbox{ , } \int_{\Omega} 
\mid \frac{D\eu}{\mu} \mid^{p} \mbox{ d}\mu < + \infty \right\}$$
denotes the Sobolev space of order $p$ with respect to the measure $\mu$. Obviously 
$$W^{1,p}_{\leb} (\Omega, R^{k}) = W^{1,p} (\Omega, R^{k})$$
Also, for any countably rectifiable surface $K = \overline{K}$ in $\Omega$ let consider the 
measure $\mu = \leb + \hen_{|_{K}}$; then $W^{1,p}_{\mu} (\Omega, R^{k})$ is a subset of 
$\eSBV(\Omega, R^{k})$ and the inclusion is continuous.Therefore the following formula holds:
\begin{equation}
\langle D\eu , \phi \rangle = \int_{\Omega} \nabla \eu \cdot \phi \mbox{ d}x + \int_{K} [\eu] 
\otimes \en \cdot \phi \mbox{ d}\hen \mbox{ \ \ } \forall \phi \in C^{\infty}(\Omega , 
R^{n \times k})
\label{stokes}
\end{equation}
Moreover 
$$W^{1,p}_{\mu} (\Omega, R^{k}) \subset W^{1,p} (\Omega \setminus K, R^{k})$$
In [DGCL] is proved that 
$$W^{1,2}_{\mu}(\Omega, R^{n}) \cap L^{\infty}(\Omega , R^{n}) = 
W^{1,2}(\Omega \setminus K , R^{n}) \cap 
L^{\infty}(\Omega , R^{n})$$

$W^{1,p}_{\mu} (\Omega, R^{k})$ is a reflexive Banach space endowed with the norm 
$$ \| \eu \|_{W^{1,p}_{\mu}} = \| \eu \|_{L^{p}_{\mu}} + \| \frac{D\eu}{\mu} \|_{L^{p}_{\mu}}$$
The weak convergence, denoted by $\eu_{h} \rightharpoonup \eu$, is equivalent to: 
$$ \eu_{h}  \rightharpoonup \eu \mbox{ in } L^{p}_{\mu} \mbox { , } \frac{D\eu_{h}}{\mu} \weak  
\frac{D\eu}{\mu} \mbox{ in } L^{p}_{\mu}$$
Another weak convergence, denoted by $\eu_{h} \weakdown \eu$, may be defined by:
$$ \eu_{h} \weakdown \eu \mbox{  } 
\Longleftrightarrow \left\{ \begin{array}{l}
                            \eu_{h} \rightarrow \eu \mbox{ in } L^{1}(\Omega, R^{k}) \\
                            \frac{D\eu_{h}}{\mu} \weak  \frac{D\eu}{\mu} \mbox{ in }
 L^{p}_{\mu}(\Omega, R^{k}
                            \end{array} \right. $$

Let consider the decomposition $\mu = a \leb + \mu^{s}$ and the functional 
$$F(x) = \int_{\Omega} f(\frac{\nabla \eu}{a} ) a \mbox{ d}x$$
where $f$ is a quasiconvex function, i.e.
\begin{equation}
f(z) \leq \int_{[0,1]^{n}} f(z + D\phi (y) ) \mbox{ d}y \mbox{ , \ \ \ } \forall z 
\in R^{n \times k} \mbox{ and } \forall \phi \in W_{0}^{1,\infty} (\Omega,R^{k})
\label{quasiconv}
\end{equation}

The following result has been proved in [AFB]:

\vspace{.5cm}

{\bf Theorem 2.1.} ( Ambrosio, Buttazzo, Fonseca) :{\it Let $p > 1$ and $f$ a quasiconvex 
function such 
that} 
$$ 0 \leq f(z) \leq C(1+ \mid z \mid^{p}) \mbox{ \ \ \ } \forall z \in R^{n \times k}$$ 
{\it If $\frac{D\mu}{dx} \in L^{\infty}$ then the functional $F$ is  seq.  lower semicontinuous on 
$W^{1,p}_{\mu}(\Omega, R^{k})$ with respect to the convergence $\weakdown$.}

\vspace{.5cm}

If $f$ is a convex function then $F$ is convex too. In the case $\mu = \leb + \hen_{|_{K}}$, if 
$\hen ( K) < +\infty$ then for any sequence $(\eu_{h})_{h} \subset W^{1,p}_{\mu}(\Omega, R^{k})$
bounded in the $L^{\infty}$ norm the convergence $\weakdown$ implies the convergence $\weak$. 
Therefore $F$ is  seq.  lower semicontinuous on every set 
$$W^{1,p}_{\mu}(\Omega, R^{k},M) = \left\{ \eu \in W^{1,p}_{\mu}(\Omega, R^{k}) \mbox{ : } \| 
\eu \|_{\infty} \leq M \right\}$$ 
This is sufficient to prove the existence of a minimizer of $F$  over the subspace 
$W^{1,2}_{\mu} (\Omega, R^{k}, \eu_{0})$ of all $\eu$ with the prescribed value $\eu_{0}$ on 
$\paromega$ ( $\eu_{0}$ in $L^{\infty}$).
The coercivity of $F$ is equivalent to:
$$ \lim_{\| \eu \|_{W^{1,2}_{\mu}(\Omega, R^{k}, \eu_{0})} \rightarrow + \infty} F(\eu) = + 
\infty$$
In the simplest case, corresponding to the Dirichlet problem for a fractured elastic body 
with natural configuration $\Omega \setminus K$, $n=k$, $f$ has the form:
$$ f(z) = \frac{1}{2} \eC z : z $$ and the property 
$$ f(z) \geq c \mid z^{sim} \mid^{2} \mbox{ \ \ \ } \forall z \in R^{n \times n}$$ where $\eC$ 
is the linear elasticity 4-order tensor, i.e. it has the properties:
$$\eC_{ijkl} = \eC_{jikl} = \eC_{ijlk} = \eC_{klij}$$
 In this case (and with the hypotheses on $\mu$ and $K$) 
the following inequality holds:
$$ \forall \eu \in W^{1,p}_{\mu}(\Omega, R^{k},M) \mbox{ \ \ } \| \eu \|_{W^{1,2}} \leq \| 
\eu \|_{W^{1,2}_{\mu}} \leq \| \eu \|_{W^{1,2}} + C(M,K)$$
So $F$ is coercive whenever the same $F$ defined over $W^{1,2} (\Omega \setminus K, R^{k})$ is 
coercive. The conclusion is:

\vspace{.5cm}

{\bf Theorem 2.2.} : {\it Let $\Omega$ be an open bounded subset of $R^{n}$, with Lipschitz 
boundary and 
$K = \overline{K} \subset \Omega$ a countably rectifiable hypersurface with $\hen ( K) < +\infty$. 
Let $\mu = \leb + \hen_{|_{K}}$ and $\eu_{0} \in  W^{1,2}_{\mu}(\Omega, R^{n}) \cap 
L^{\infty}(\Omega , R^{n})$. Then the functional:
$$F: W^{1,2}_{\mu} (\Omega, R^{k}, \eu_{0}) \rightarrow R$$
$$F(\eu) = \frac{1}{2} \int_{\Omega} \eC \nabla \eu : \nabla \eu \mbox{ d}x \mbox{ \ \ ,}$$ 
where $\eC$ is a symmetric positive definite 4-order elasticity tensor,
has a minimizer, unique to an element of 
$$R^{1,2}_{\mu}(\Omega, R^{n},0) = \left\{ \eu \in W^{1,2}_{\mu}(\Omega, R^{k}, 0) 
\mbox{ : } \epsilon(\eu) = 0 \right\}$$
Any minimizer $\eu$ has the properties:
\begin{itemize}
\item $\eu \in W^{2,p}_{loc}(\Omega \setminus K, R^{n})$ $\forall p \geq 2$
\item $div \eC \nabla \eu = 0$ in $\Omega$ in the sense of distributions
\item $\left( \eC \nabla \eu \right) \en = 0$ on $K$
\item $\forall \phi \in C^{\infty}(\Omega,R^{n \times n})$ 
$$\int_{\Omega} div \phi \mbox{ }\cdot \eu \mbox{ d}x + \int_{\Omega} \nabla \eu : \phi \mbox{ 
d}x + \int_{K} \left[ \eu \right] \cdot \left( \phi \en \right) \mbox{ d}\hen = 
\int_{\paromega} \left( \phi \en \right) \cdot \eu_{0}  \mbox{ d}\hen$$
\end{itemize} }

\vspace{.5cm}

The uniqueness of the minimizer $\eu$ follows from the expression of $F$ in a classical way. The 
properties of $\eu$ come from the Euler equation:
$$\int_{\Omega} \eC \nabla \eu : \nabla \ev \mbox{ d}x = 0 \mbox{ \ \  \ } \forall \ev \in 
W^{1,2}_{\mu} (\Omega,R^{n},0)$$ and from the remarks from the beginning of the section. 

The kinematic of a fractured elastic medium is therefore well settled in the frame of the Sobolev 
space $W^{1,2}_{\mu} (\Omega,R^{n})$. The minimizer $\eu$ is a physical solution of the Dirichlet 
problem since the formula (\ref{stokes}) shows that there are no forces or moments concentrated 
at the edge of the crack $K$ (the border of $K$).

The space of admissible stresses is (according to [DP]) :
$$Y_{K}(\Omega) = \left\{ \sigma = \sigma^{T} \in L^{2}(\Omega,R^{n \times n}) \mbox{ : } 
div \sigma = 0 \mbox{ , } \sigma \en = 0 \mbox{ on } K \right\}$$
Here $\sigma \en$ means the trace of $\sigma$ on $K$ multiplied by the normal $\en$. 
The definition is good because of the following proposition:

\vspace{.5cm}

{\bf Proposition 2.3.} (Del Piero) : {\it Let $Y(\Omega)$ be the space:
$$Y(\Omega) = \left\{ \sigma = \sigma^{T} \in L^{2}(\Omega,R^{n \times n}) \mbox{ : } 
div \sigma \in L^{2}(\Omega,R^{n}) \mbox{ , } \right\}$$ and $K$ any smooth surface in $\Omega$. 
Denote by $\left[ \sigma \en \right]$ the jump of $\sigma \en$ over $K$. Then 
$\left[ \sigma \en \right] = 0$.}

\vspace{.5cm}

The machinery of convex duality (see [M]) holds in this case. Let consider the following separate 
bilinear form:
$$ \langle , \rangle : L^{2}(\Omega,R^{n \times n}_{sym} \times L^{2}(\Omega,R^{n \times n}_{sym}
 \rightarrow R$$
$$ \langle \sigma , \epsilon \rangle = \int_{\Omega} \sigma : \epsilon \mbox{ d}x$$
The polar of the convex functional
$$W:L^{2}(\Omega,R^{n \times n}_{sym} \rightarrow R \mbox{ \ \ \ } W(\epsilon) = 
\frac{1}{2} \int_{\Omega} \eC \epsilon : 
\epsilon \mbox{ d}x$$ is defined by
$$W^{*} : L^{2}(\Omega,R^{n \times n}_{sym} \rightarrow R \cup \left\{ +\infty \right\}$$
$$W^{*}(\sigma) = \sup \left\{ \langle \sigma , \epsilon \rangle - W(\epsilon) \mbox{ : } 
\epsilon \in L^{2}(\Omega,R^{n \times n}_{sym} \right\} = 
\frac{1}{2} \int_{\Omega} \eC^{-1} \sigma : \sigma \mbox{ d}x$$
The following inequality holds:
$$W^{*}(\sigma') + W(\epsilon") \geq \langle \sigma',\epsilon" \rangle$$ and the equality 
is attained whenever $\sigma' = \eC \epsilon"$. If $\eu$ denotes a minimizer of $F$, then 
$\eC \nabla \eu \in Y_{K}(\Omega)$.

For any $\sigma \in Y_{K}(\Omega)$ and for any $\ev \in W^{1,2}_{\mu}(\Omega, R^{n}) 
\cap W^{2,p} (\Omega \setminus K, R^{n})$
$$\langle \sigma , \epsilon(\ev) \rangle = \int_{\Omega} \sigma_{ij} \ev_{i,j} \mbox{ d}x + 
\int_{K} [\eu_{i}] (\sigma \en )_{i} \mbox{ d}\hen = \left( D\eu , \sigma \right) =$$  $$=
 - \int_{\Omega} \eu_{i} ( div \sigma )_{i} \mbox{ d}x + \int_{\paromega} (\sigma \en) \cdot 
\eu_{0} \mbox{ d}\hen = \int_{\paromega} (\sigma \en) \cdot 
\eu_{0} \mbox{ d}\hen $$ This, together with the inequality involving $W$ and $W^{*}$, give 
the following

\vspace{.5cm}

{\bf Proposition 2.4.} : {\it For any $\sigma \in Y_{K}(\Omega)$ and for any minimizer of $F$ the 
following inequality holds: 
$$\int_{\paromega} (\sigma \en) \cdot \eu_{0} \mbox{ d}\hen - \frac{1}{2} \int_{\Omega} 
\eC^{-1} \sigma : \sigma \mbox{ d}x \leq F(\eu)$$ }

\section{Introduction}
\indent

Any Mumford-Shah functional   is related to a perturbed area functional. 
The result of the first section shows that, in 
particular cases which model the fiber-matrix decohesion in some composite materials, 
the discontinuity set of the minimizer of the classical Mumford-Shah functional 
minimizes the area functional too.

The decohesion between the fibers and the matrix in a composite material can be modeled as fracture 
of the matrix. In the case of a composite with rigid fibers and elastic matrix the fracture is 
brittle and it has the peculiar property that a crack may appear on the boundary of the matrix. If 
this happens then the matrix no longer satisfies the boundary displacement imposed by the 
fibers and the decohesion appears. The problem is therefore to model the appearance of a crack 
in an elastic body submitted to a given boundary displacement. The purpose of the first section 
is to 
show that in the case of antiplane displacement the energetic approach to the problem can model 
the decohesion, i.e. the crack predicted by the model may lie on the boundary of the body. 
We prove that in a particular class of geometries and boundary conditions the only crack that may 
appear is a geodesic (i.e. length minimizing set). The result is quantitative and experiments can
be made in order to validate the energetic model.

In the second section we consider two phenomena that occur in the damage of an 
elastic body: brittle (or quasi-brittle) damage and brittle fracture. These two damage 
mechanisms are the macroscopic effects of the appearance of a crack. That crack may be 
formed by a macro-crack and a large number of meso-cracks; the existence of the latter explains 
the weakening of the elastic macro-properties of the body. Brittle damage is modeled as 
something occurring in the volume of the body and the theory predicts the existence of fine 
mixtures between the  damaged and sound states (see [FMu] and [FMa]). This is a mathematical 
consequence of the fact that for a minimizing sequence $(\eu_{h},A_{h})$ of the total energy 
functional
\begin{equation}
E(\ev,A) = \int_{\Omega} w(\nabla \eu, A) \mbox{ d}x + \gamma \int_{\Omega} \chi_{A} \mbox{ d}x
\label{damagef}
\end{equation}
(where $A$ is the damaged region) the sequence $(\chi_{A_{h}})$ might not converge to a 
characteristic function of a set, but to an $f:\Omega \rightarrow [0,1]$. The reason is that there 
is no control of the perimeter of the damage region. 

The energetic model of 
crack appearance (in mode 3, antiplane) predicts a priori non-smooth fractures. 
There are three fundamental facts about the cracks $K$ predicted by this model: 
\begin{itemize}
\item  1. $K$ is not a curve that fills a volume: $\hen(\overline{K} \setminus K) = 0$
\item 2. there is a piecewise smooth manifold containing $K \setminus S$ where $S$ is a closed  
$\hen$-negligible set (see [DS], [AP] and [AFP]).
\item 3. the curvature of $K$ is proportional to the jump of the bulk energy density $[w]$ (this 
relation is true only in a formal sense due to the lack of smoothness of $K$) 
\end{itemize}
Our purpose is to couple brittle (brutal and total) damage with brittle fracture. This can be done 
in  two ways; the first is to add to the crack $K$ the essential frontier of the damaged zone
 and add the  volume of the damaged zone in the expression of the energy functional. 
The surface term in the 
expression of the energy will  control then the perimeter of the damaged zone as well as the length
 of the  crack. In this model no fine mixtures appear.

 The second way to couple brittle damage with 
brittle fracture is to renounce to control the perimeter of the damage zone; then, in the 
portions of the crack embedded in sound regions (because now fine mixtures may appear), the 
curvature of the crack has an upper bound. This is due to the fact that  the volume of the damaged
 zone  controls the jump of the bulk density hence the 
curvature of the crack. In this model we choose to extend the class of $n-1$ smooth manifolds to 
the class of integral varifolds.

\section{Decohesion in composites and Mumford-Shah functionals}
\indent

\vspace{.5cm}

The idea of the energetic models of crack appearance is to enlarge the class of admissible 
displacements of the body from the space $W^{1,2}(\Omega,R^{n})$ to $\eSBV(\Omega,R^{n})$ or even 
to $\eSBD(\Omega,R^{n})$. The jump of the displacement (let say from $\eSBV(\Omega,R^{n})$) 
represents a crack; in this way all the displacements corresponding to some arbitrary crack in the 
body lie  in $\eSBV(\Omega,R^{n})$. In the preliminaries we see that the correct setting of 
the equilibrium Dirichlet problem for an elastic cracked body is in $W^{1,2}_{\mu} (\Omega, 
R^{n})$, with $\mu = \leb + \hen_{|_{K}}$ ($K$ is the crack); this is a subspace of 
$\eSBV(\Omega,R^{n})$. Suppose that on a non-negligible part $\gau$ of the boundary of the body
$$\paromega = \overline{\gau} \cup \overline{\gaf} \mbox{ , } \gau \cap \gaf = \emptyset 
\mbox{ , } 0 < \hen(\gau) \cdot \hen(\gaf) < +\infty$$ the boundary displacement $\eu_{0}$ is 
imposed and $\gaf$ is stress free. For technical reasons $$\eu_{0} \in \eSBV_{loc}(R^{n},R^{n}) \cap 
L^{\infty}(R^{n},R^{n}) \mbox{ \ \ , \ \ } \hen(\gaf \setminus \eS_{\eu_{0}} ) = 0$$ Here $\eS_{\eu}$
 denotes 
the jump set of $\eu$; because $\eu \in \eSBV_{loc}(R^{n},R^{n})$ the jump set which is usually 
denoted by $\eJ_{\eu}$ coincides modulo a $\hen$-negligible set with the complement of the  
Lebesgue set $\eS_{\eu}$. The space of admissible displacements is
$$\eSBV(\Omega, R^{n}, \eu_{0}) = \left\{ \eu \in \eSBV_{loc}(R^{n},R^{n}) \mbox{ : } 
\eu = \eu_{0} \mbox{ in } R^{n} \setminus \overline{\Omega} \right\}$$
We are searching for a minimizer $\eu \in \eSBV(\Omega, R^{n}, \eu_{0})$ of the functional $I$ 
over the same space, where the total energy $I$ is a Mumford-Shah type functional:
$$I(\ev) = \frac{1}{2} \int_{\Omega} \eC \nabla \ev : \nabla \ev \mbox{ d}x + G \hen(\eS_{\ev})$$

We remark that $I$ is well defined over $\eSBD(\Omega)$ if one replaces $\eS_{\eu}$ with 
$\eJ_{\eu}$. The existence of a minimizer of $I$ in this 
space is proved in [BCDM].  

In the case of antiplane displacements --- $n=3$ and $\eu(x_{1}, x_{2}, x_{3}) = (0, 
0,\eu(x_{1},x_{2}))$ --- the correct setting of the problem is to consider $\Omega \in R^{2}$ and 
$\eu \in \eSBV(\Omega, R, \eu_{0})$. The expression of the total energy $I$ becomes
$$I(\ev) = \frac{1}{2} \int_{\Omega} \mid \nabla \ev \mid^{2} \mbox{ d}x + G \hen(\eS_{\ev})$$
It is known from [Amb] and [DGCL] that the latter $I$ is $L^{1}_{loc}$ lower compact so there is 
at least a minimizer $\eu$ of $I$ in the space $\eSBV(\Omega, R, \eu_{0})$. Much is known 
about the geometric properties of the crack predicted by this model, i.e. $\eS_{\eu}$; among them 
we mention the fact that $\eS_{\eu}$ is a subset of a smooth hypersurface, excepting a 
$\mathcal{ H}^{1}$- negligible set (recently proved in [DS]). 
That push us to return to the original Mumford-Shah functional and minimization problem. 
Consider the functional 
$$J(\ev,K) = \frac{1}{2} \int_{\Omega} \eC \nabla \ev : \nabla \ev \mbox{ d}x + G \hen(K 
\setminus \gaf)$$ defined 
over the set $M(\eu_{0})$ of all $(\ev,K)$ such that $K = \overline{K} \subset \overline{\Omega}$ is 
a piecewise $C^{1}$ 
hypersurface, $\gaf \subset K$ and $\ev \in W^{1,2}_{\mu}(\Omega,R^{n})$, 
where $\mu = \leb + \hen_{|_{K}}$. Also 
$\ev = \eu_{0}$ on $\gau \setminus K$. With the help of $J(\cdot,\cdot)$  we can  define a 
perturbed area functional over all piecewise $C^{1}$ closed hypersurfaces:
$$J(K) = \min \left\{ J(\ev,K) \mbox{ : } (\ev,K) \in M \right\}$$
$J(\cdot)$ may or may not have a minimizer and this is the reason for using the functional $I$ first 
proposed by Ambrosio. However it seems natural to us that for a big $G$ the minimizer (if it exists)
 of $J(\cdot)$ minimizes also the area $\hen( \cdot )$. Such a situation is described in the sequel.

\vspace{1.5cm}

For any smooth $K$ (i.e. with the properties listed several times before) denote by $\eu_{K}$ the 
displacement of the fractured body $\Omega \setminus K$ under the imposed displacement $\eu_{0}$.
$$J(\eu_{K}, K) = J(K)$$
We restrict to the case of antiplane displacements so $\eu$ is a scalar function and the bulk 
energy is a quadratic expression in $\nabla \eu$. Suppose that $\eu_{\gaf} \in C^{1}$ and 
$0 < c < \mid \nabla \eu_{\gaf} \mid < C$ everywhere in $\Omega$. Define then  $\ueu$ to be 
any function with the properties:
$$\frac{\partial \ueu}{\partial x_{1}} = \frac{1}{\| \eu_{0} \|_{\infty}} \frac{\partial 
\eu_{\gaf}}{\partial x_{2}}$$
$$\frac{\partial \ueu}{\partial x_{2}} = - \frac{1}{\| \eu_{0} \|_{\infty}} \frac{\partial 
\eu_{\gaf}}{\partial x_{1}}$$
The level surfaces of $\ueu$ form a congruence of curves in $\Omega$. The free surface $\gaf$ is 
tangent to the congruence. This congruence defines a system of open (smooth) neighbourhoods 
$V(\ueu)$ by:
$$\forall A \in V(\ueu) \mbox{ \ \ \ }  \partial A \  \setminus \gau \mbox{ is locally a level 
set of } \ueu$$
For any $A \in V(\ueu)$ we will denote by $\partial_{\eu} A$ the part of the boundary of $A$ made 
by curves $\eu = ct$.
Let $K$ be a piecewise smooth curve and $\Omega' \in V(\ueu)$ such that $K \setminus \gaf 
\subset \Omega'  \cup 
\partial_{\eu} \Omega^{'}$. 
Recall that 
$$J(K) =  \frac{1}{2} \int_{\Omega} \mid \nabla \eu_{K} \mid^{2} \mbox{ d}x + G \mathcal{ H}^{1}(K)$$ 
The following field is an admissible stress (i.e. it belongs to $Y_{K}(\eu_{\gaf})$, see the 
preliminaries):
$$ \sigma = \left\{ \begin{array}{ll}
                    \nabla \eu_{\gaf} & \Omega \setminus \Omega' \\
                       0      & \Omega'
                    \end{array} \right. $$
Proposition 3.4 (which affirms the familiar minimum principle in stress) gives  
$$\frac{1}{2} \int_{\Omega} \mid \nabla \eu_{K} \mid^{2} \mbox{ d}x \geq \int_{\gau} (\sigma \en) 
\cdot \eu_{0} \mbox{ d}\mathcal{ H}^{1} - \frac{1}{2} \int_{\Omega} \mid \sigma \mid^{2} \mbox{ d}x$$
$$J(\gaf) - J(K) \leq - \frac{1}{2}\int_{\Omega'} \mid \nabla \eu_{\gaf} \mid^{2} + 
\int_{\gau \cap \overline{\Omega'}} (\nabla \eu_{\gaf} \en) \cdot \eu_{0} - G \mathcal{ H}^{1}(K 
\setminus \gaf)$$
Denote by $\tau$ the clockwise tangent to $\paromega'$. After some work we find
$$J(\gaf) - J(K) \leq \frac{1}{2} \| \eu_{0} \|_{\infty} \int_{\paromega'} \eu_{\gaf} \cdot (
\nabla \ueu \tau ) - G \mathcal{ H}^{1}(K 
\setminus \gaf)$$ 
Suppose now that $\gau$ has two connected components $\gau^{1}$ and $\gau^{2}$ and $\eu_{0}$ is 
taken like this:
$$\eu_{0} = 0 \mbox{ on } \gau^{1} \mbox{ and } \eu_{0} = \Delta = \| \eu_{0} \|_{\infty} 
\mbox{ on } 
\gau^{2}$$

Consider now the projection, called $P$, along the congruence $\ueu = ct.$ on $\gau^{2}$. Then
$$\frac{1}{2} \int_{\partial_{\eu} \Omega^{'}}  (\nabla \eu_{\gaf} \en) \cdot \eu \mbox{ d}
\mathcal{ H}^{1} \geq \frac{1}{2} \Delta^{2} \ V \ueu (P(K))$$
 where $V \ueu (M)$ is the variation of $\ueu$ over the set $M$, $\mathcal{ H}^{1} (M) < + \infty$. 
Moreover
$$\inf \left\{ \frac{1}{2} \int_{\partial_{\eu} \Omega^{'}}  (\nabla \eu_{\gaf} \en) \cdot \eu 
\mbox{ d}\mathcal{ H}^{1} \mbox{ : } \Omega^{'} \in V(\ueu), K \setminus \gaf \subset \Omega^{'}
\right\} = 
\frac{1}{2} \Delta^{2} \ V \ueu (P(K))$$ Remark also that 
$$V\ueu(K \setminus \gaf) \geq V\ueu(P(K))$$ Take $\overline{D}(K)$ to be the intersection of all
 members of 
$V(\ueu)$ containing $K \setminus \gaf$ and denote by $\Gamma(K)$ any of the smallest geodesics in  
$\overline{D}(K)$ which separate the two connected components of $\partial \overline{D}(K) \cap 
\gau$. Then
$$ \mathcal{ H}^{1}(\Gamma(K)) \leq \mathcal{ H}^{1}(P(K)) \mbox{ , } \mathcal{ H}^{1}(\Gamma(K)) \leq 
\mathcal{ H}^{1}(K \setminus \gaf) \mbox{ , } V\ueu (P(K)) \leq V\ueu (\Gamma(K))$$
The following inequality holds because on the assumption on $\mid \nabla \eu_{\gaf} \mid$:
$$\frac{1}{2} \Delta^{2} V\ueu (K) - G\mathcal{ H}^{1}(K \setminus \gaf) \leq 
\mathcal{ H}^{1}(K \setminus \gaf)\left( \frac{1}{2} \Delta^{2} C - G \right)$$
Therefore if $$\Delta^{2} \leq \frac{2G}{C}$$ then 
$$J(\gaf) - J(K) \leq \frac{1}{2} \Delta^{2} V\ueu (K) - G\mathcal{ H}^{1}(K \setminus \gaf) \leq 0$$

Also, remark that
$$\sup \left\{ \frac{1}{2} V \ueu (P(K)) - G\mathcal{ H}^{1}(K \setminus \gaf) \mbox{ : } 
K = \overline{K} \right\} = \sup \left\{ \frac{1}{2} V \ueu (\Gamma(K)) - G\mathcal{ H}^{1}
(\Gamma(K)) \mbox{ : } K = \overline{K} \right\}$$
$$\frac{1}{2} \Delta^{2} V\ueu (\Gamma(K)) - G\mathcal{ H}^{1}(\Gamma(K)) \geq 
V\ueu(\Gamma(K)) \left( \frac{1}{2} \Delta^{2} - G\frac{1}{c} \right)$$
Therefore if 
$$\Delta^{2} \geq \frac{2G}{c}$$
then the following measure over the $\sigma$-algebra over $V(\ueu)$ is a positive one:
$$\Omega^{'} \mapsto \frac{1}{2} \int_{\Omega^{'}} \mid \nabla \eu_{\gaf} \mid^{2} \mbox{ d}x - 
G \mathcal{ H}^{1}(\Gamma(\partial_{\eu} \Omega^{'}))$$
and
$$J(\gaf) - J(K)  \leq \frac{1}{2} \int_{\overline{D}(K)} \mid \nabla \eu_{\gaf} \mid^{2} 
\mbox{ d}x - G \mathcal{ H}^{1}(\Gamma(\partial_{\eu} \overline{D}(K)))$$
The maximum of the right term from the previous inequality is attained in $\overline{D}(K) = \Omega$
and for any smallest geodesic
in $\Omega \cup \partial_{\eu} \Omega$ separating the two components of $\gau$ this inequality 
becomes an equality.

\vspace{.5cm}

{\bf Proposition 1} : {\it  With the assumptions on $\Omega$ and $\eu_{0}$ mentioned before, 
there exist two positive constants $m \leq M$, depending only on $\Omega$, $\gau$, $G$, such that: 
\begin{itemize}
\item if $$\| \eu_{0} \|_{\infty}^{2} \ \leq \ m$$ then $J(\gaf) \leq J(K)$ for all admissible 
$K$; this means that no crack appears under the imposed antiplane displacement $\eu_{0}$,
\item if  $$\| \eu_{0} \|_{\infty}^{2} \ \geq \ M$$ then  any smallest geodesic in 
$\overline{\Omega}$ which separates $\gau^{1}$ from $\gau^{2}$ is a minimizer of $J$.
\end{itemize}
} 

\vspace{.5cm}

The constants $m$ and $M$ can be precisely determined in some particular cases. For example 
if $\Omega = B(0,R) \setminus \overline{B(0,r)}$ is a ring and $\eu_{0} = 0$ on 
$\partial  B(0,R)$ and $\eu_{0} = const.$ on $\partial  B(0,r)$ then 
$$m = M = r \ \ ln \frac{R}{r}$$ The only crack that may appear in this case 
is $\partial  B(0,r)$. 

If $\Omega = (0,a) \times (0,L)$ is a rectangle with the free surfaces 
$\left\{ 0 \right\} \times (0,L)$, $\left\{ a \right\} \times (0,L)$ and $\eu_{0} = 0$ on 
$(0,a) \times \left\{ 0 \right\}$, $\eu_{0} = const.$ on  $(0.a) \times \left\{ L \right\}$, 
then (calculus done by J.-J. Marigo) 
$$m = M = L$$ The only crack that may appear in this case is any line 
$(0,a) \times \left\{ x \right\}$, with $x \in [0,L]$.

 There are two types of mechanical experiments that can be made 
in these cases in order to validate the model. One may try first to compare, for a given brittle 
material and geometry, the Griffith material constant with the constant $G$ determined from the 
critical displacement measurements. Second, for the same material but for different geometries, 
the experimental values of $G$ can be compared in order to see if $G$ is a material constant.

\vspace{1.5cm}

\section{Perturbed area functionals in brittle damage mechanics}
\indent

Since the Mumford-Shah functional seems to be intimately connected with the area functional then 
the methods and objects used in the variational problems involving the area functional may be 
useful in some free boundary problems. We present in this section two different models of 
brittle damage coupled with brittle fracture. In the first we use the $\eSBD$ approach but 
in the second we enlarge the class of  $n-1$ smooth  manifolds to the class of integral 
varifolds.  This is not a new idea; in the [AFP] and [AP] the authors use varifold techniques 
to prove the regularity of the minimizer of the Mumford-Shah functional.  

\vspace{.5cm}

Let consider the following space:
$$\mathcal{ M} = \left\{ (\eu,A) \mbox{ : } \mathcal{ P}(A,\Omega) < + \infty \mbox{ , } \eu \in 
\eSBD(\Omega) \mbox{ , } E \left( \eu \chi_{A^{c}} \right) \ll E \chi_{A} \right\}$$
We define on $\mathcal{ M}$ the  energy functional ($A$ is the sound region; from here the minus sign in
the expression of the energy):
$$J(\eu,A) = \int_{A} \left[ w(\epsilon (\eu)) - \gamma \right] \mbox{ d}x + G \hen \left( 
\eJ_{\eu} \cup \partial^{*} A \right)$$

\vspace{.5cm}

{\bf Theorem 1} {\it Let $(\eu_{h},A_{h}) \subset \mathcal{ M}$ be a sequence such that for all $h \in $ 
{\bf N}} 
\begin{equation}
\left\{ \| \eu_{h} \|_{\infty} + J(\eu_{h},A_{h}) \right\} \leq C
\label{th1}
\end{equation}
{\it Then one can find a subsequence $(\eu_{h_{k}}, A_{h_{k}})$ and an element $(\eu,A) \in 
\mathcal{ M}$ such that:}
$$\eu_{h_{k}} \rightarrow \eu \mbox{ strong in } L^{1}_{loc}(\Omega,R^{n})$$
$$\chi_{A_{h_{k}}} \rightarrow \chi_{A} \mbox{ in } \eBV(\Omega,R)$$
$$ J(\eu,A) \leq \liminf_{k \rightarrow \infty} J(\eu_{h_{k}}, A_{h_{k}})$$

\vspace{.5cm}

Before giving the proof,  remark that this theorem gives the existence of the minimizers 
of $J$ in the set $$\mathcal{ M}(\eu_{0}) = \left\{ (\eu,A) \mbox{ : } (\eu_{|_{\Omega}},A) \in 
\mathcal{ M} \mbox{ , } \eu \chi_{A}  = \eu_{0} \mbox{ on } R^{n} \setminus \overline{\Omega} \right\}$$
where $\eu_{0} \in \eSBD_{loc}(R^{n}) \cap L^{\infty}(R^{n},R^{n})$. Remark also that for any 
minimizing pair $(\eu,A)$ $\| \eu \chi_{A} \|_{\infty} \leq \| \eu_{0} \|_{\infty}$. By a 
translation argument one can prove then that the functional
$$K(\eu,A) = \int_{\Omega} \left[ w(\epsilon(\eu \chi_{A} ) ) - \gamma \chi_{A} \right] \mbox{ d}x 
+ G \hen(\eJ_{\eu \chi_{A}})$$ has minimizers in $\mathcal{ M}(\eu_{0})$.

\vspace{.5cm}

Proof of the Theorem 1. By the compactness theorem for sets with finite perimeter it follows that 
we can extract a subsequence of $(\eu_{h},A_{h})$ (denoted by the same name as the initial 
sequence) such that 
$$\chi_{A_{h}} \rightarrow \chi_{A}$$ in $\eBV(\Omega,R)$. Since $(\eu_{h},A_{h}) \in \mathcal{ M}$ 
then $\eu$ is a linear antisymmetric function out of $A_{h}$.
 The assumptions on $w$, (\ref{th1}) and $\leb(\Omega) < + \infty$ imply that for all $h$:
$$\int_{\Omega} \mid \eu_{h} \mid \mbox{ d}x + \mid E^{j} \eu_{h} \mid (\Omega) + c 
\int_{\Omega} \mid \epsilon ( \eu_{h} ) \mid^{2} \mbox{ d}x + \hen(\eJ_{\eu_{h}}) \leq C^{'}$$
By the compactness theorem proved in [BCDM] we can extract a subsequence such that (keeping the 
notation):
$$\eu_{h} \rightarrow \eu \mbox{ strong in } L^{1}_{loc}(\Omega, R^{n})$$
$$E^{j}\eu_{h} \rightarrow E^{j} \eu \mbox{ weak * as measures}$$
$$\epsilon(\eu_{h}) \rightarrow \epsilon(\eu) \mbox{ weak in } L^{1}(\Omega, R^{n \times 
n}_{sym})$$
Therefore $E\eu_{h} \rightarrow E\eu$ weak * as measures and $\chi_{A_{h}^{c}} \rightarrow 
\chi_{A^{c}}$ in $eBV(\Omega,R)$, hence
$$E\left( \eu_{h} \chi_{A_{h}^{c}} \right) \rightarrow E \left( \eu \chi_{A^{c}} \right) 
\mbox{ weak * as 
measures}$$
(\ref{th1}) implies again that the densities $\frac{E\left( \eu_{h} \chi_{A_{h}^{c}} \right)}
{E \chi_{A_{h}}}$ are uniformly bounded in $L^{\infty}$, therefore one of the Hutchinson theorems 
(see [Hut], 4.4.2.)  
assures us that $$E \left( \eu \chi_{A^{c}} \right) \ll E \chi_{A}$$ 
The lower semicontinuity of $J$ is a straightforward consequence of the previous facts.

\vspace{1.cm}

Recall that a $(n-1)$  varifold is a positive Radon measure on the bundle 
$\Omega \times G(n,n-1)$  
with 
the grassmanian $$G(n,n-1) = \left\{ \eS= \left( \eS_{ij} \right) \mbox{ : } \eS_{ij} = \delta_{ij} 
- 
\en_{i} \en_{j} \mbox{ \ ,\ } \en \in S^{n} \right\}$$ as fiber in any $x \in \Omega$. The weight 
of the varifold $\eV$ is $$\| \eV \| (A) = \eV \left( \left\{ (x,\eS) \mbox{ : } x \in A \right\} 
\right) \mbox{ \ \ \ } \forall A \in \eB(\Omega)$$ 
$\eV^{(x)}(\eS)$ denotes the density of $\eV$ in the fiber of $x \in \Omega$ of the mentioned 
bundle. The class of unoriented $n-1$ integer varifolds is denoted by $IV(\Omega,n-1)$. One can 
consider also the class of oriented $n-1$ varifolds (Radon measures over the bundle 
$\Omega \times G^{o}(n,n-1)$); if $q$ is the standard map from $ G^{o}(n,n-1)$ to $G(n,n-1)$ which 
associates to any oriented hyperplane the unoriented one, then $q^{\#}$ maps oriented varifolds 
onto unoriented ones. The class of integer oriented $n-1$ varifolds is denoted by 
$IV^{o}(\Omega,n-1)$.

The first variation of the varifold $\eV$ is denoted by $\delta \eV$ (with the total variation as 
measure 
$\mid \delta \eV \mid$) and it has the 
definition: for any $g \in C^{1}_{0}(R^{n},R^{n})$
$$\delta \eV (g)  = \int g_{i,j}(x) \eS_{ij} \mbox{ d}\eV(x,S)$$

\vspace{.5cm}

In the following $(\eu,A,K)$ is an admissible triple if $\eu \in W^{1,2}_{\mu}(\Omega,R^{n})$, 
$\mu = \leb + \hen_{|_{K}}$, $K \subset A $ is a piecewise smooth surface and $\chi_{A} \in \eBV(
\Omega,R)$ open with smooth frontier. 
$A$ denotes the sound region of the body $\Omega$ (from here the minus sign in the 
further expression of the total energy). Any minimizer (if exists) of the functional $I$
$$I(\eu,A,K) = \int_{\Omega} \left( w(\nabla \eu) - \gamma \right) \chi_{A} \mbox{ d}x + 
 G\hen(K)$$ in 
the set of all admissible triples with $\eu$ satisfying a Dirichlet boundary condition
is a critical admissible 
triple in the sense that  for any 
$\eta \in W^{2,p}(\Omega,R^{n})$, $supp \  \eta$ compact in $\Omega$,
$$\delta_{\eta} I (\eu,A,K) = \liminf_{\epsilon \rightarrow 0} \frac{I(\eu.\phi_{\epsilon}^{-1}, 
\phi_{\epsilon}(A),\phi_{\epsilon}(K)) - I(\eu,A,K)}{\mid \epsilon \mid} \geq 0$$
where $\phi_{\epsilon}$ is the one-parameter group generated by $\eta$. Here $w$ is a $C^{2}$ 
strict convex positive function with the following properties:
\begin{itemize}
\item $w(B) = w(B^{sym})$ for all $B \in R^{n \times n}$,
\item there are two positive constants $c$ and $C$ such that 
$c \mid B^{sym} \mid^{2} \leq w(B) \leq C \mid B^{sym} \mid^{2}$  for all 
$B \in R^{n \times n}$.
\end{itemize}

The first variation of $I$ with respect to $\eta$ is
\begin{equation}
\delta_{\eta} I (\eu,A,K) = \int_{\Omega} \left[ \left( w(\nabla \eu) - \gamma \right) div \eta 
- \frac{\partial w}{\partial \nabla} \nabla \eu \ \nabla \eta \right] \chi_{A} \mbox{ d}x + 
G \int_{K \cup \partial A} div_{s} \eta \mbox{ d}\hen
\label{first}
\end{equation}
where $div_{s}$ is the tangent divergence with respect to a surface with normal $\en$:
$$div_{s} \eta \ = \eta_{i,i} - \eta_{i,j} \en_{i} \en_{j}$$ In fact if we see $K$ as a varifold 
then 
$$ \delta K(\eta) = \int_{K} div_{s} \eta \mbox{ d}\hen$$

Denote by $H = - div_{s} \nu$ ($\nu$ is the normal to $K$) the scalar mean curvature of $K$. 
If $(\eu,A,K)$ is a 
critical point of $I$ and if $(\eu,A,K)$ is balanced (i.e. $\sigma(\eu) = 
w_{,\nabla}(\nabla \eu) \ \in Y_{K}(A)$)
 then 
$$\left[ w(\nabla \eu ) \right] = G H $$ on $K$ in the distribution sense.  Also, a simple 
integration by parts in (\ref{first}) gives for balanced triples:
\begin{equation}
\delta_{\eta} I (\eu,A,K) = \int_{\Omega} div \left[ \left( w(\nabla \eu) - \gamma \right)  \eta 
\right] \mbox{ d}x
\label{firstb}
\end{equation}

For an $\eu_{0} \in \eSBV_{loc}(R^{n},R^{n}) \cap L^{\infty}(R^{n},R^{n})$ a triple $(\eu,A,K)$ is 
$\eu_{0}$-admissible if it satisfies the supplementary condition: 
$$\eu  \ \in \eSBV(\Omega,R^{n},\eu_{0})$$ 

We shall work with the  set $T(\eu_{0})$ of all $\eu_{0}$-admissible  triples $(\eu,A,K)$ 
satisfying the following property:
$$\exists N = N(\eu,A,K) \in R \ \    \forall \xi \in C^{0}_{c}(\Omega, [0,+\infty)) \ \ 
 N \int_{K}  \xi \mbox{ d}\hen + \int_{K} div_{s} \xi \en \mbox{ d}\hen \geq 0$$
For any element $(\eu,A,K) \in T(\eu_{0})$ $K$ has bounded scalar mean curvature, but this bound 
depends on the element. For any positive number $N$ we define the  space $T(\eu_{0},N)$ to be the 
set of all $(\eu,A,K) \in T(\eu_{0})$ such that:
$$\forall \xi \in C^{0}_{c}(\Omega, [0,+\infty)) \ \ 
 N \int_{K}  \xi \mbox{ d}\hen + \int_{K} div_{s} \xi \en \mbox{ d}\hen \geq 0$$ 
Off course
$$T(\eu_{0}) = \bigcup_{N \in R_{+}} T(\eu_{0},N)$$

For any  $(\eu,A,K) \in T(\eu_{0},N)$, if there is a subset $D$ of $\Omega$ such that 
$\int_{D\cap A} 
w(\nabla \eu) - \gamma \mbox{ d}x \geq 0$ then $$I(\eu,A,K) \geq I(\eu,A \cup D, K)$$ so the 
minimizing sequences have to be searched in the set of admissible triples with 
\begin{equation}
\int_{B_{r}(x)} w(\nabla \eu) - \gamma \mbox{ d}x \leq 0
\label{maxenerg}
\end{equation}
 for $\leb$-all $x \in \Omega$ and 
for all $r >0$ small enough. 

If we keep $(A,K)$ and vary $\eu$ we always find a minimizer of $I$ in this class. This minimizer 
is obviously balanced. Therefore we have two simple operations that lower the value of $I$. 
The question is: for a   given $N$ and a given $(\eu,A,K) \in T(\eu_{0},N)$, is there a balanced 
triple in the closure  of a $T(\eu_{0},M)$ with respect to a particular convergence, satisfying 
(\ref{maxenerg}), with lower energy $I$?

 Suppose that we repeat the "cut" and "balance" operations 
 without finding the desired triple after a finite number of iterations. We obtain a sequence 
$(\eu_{h},A_{h},K)$. For any $h \in N$ $\| \eu_{h} \|_{\infty} \leq 
\| \eu_{0} \|_{\infty}$. Therefore, up to  a subsequence, $\eu_{h} \chi_{A_{h}}$  converges weak * 
in $L^{\infty}$  to a $\eu$. 

Consider now the sequence of measure-functions (see [Hut] for definition 
and properties of measure-functions)  $(\chi_{A_{h}} \leb , \epsilon(\eu_{h}))_{h}$. It is easy 
to see that there is a constant $M$  with $$\int_{\Omega} w(\epsilon(
\eu_{h}) \  \chi_{A_{h}} \mbox{ d}x \leq M$$ for all $h \in N$. Then the 
compactness theorem 4.4.2. [Hut] implies that, up to a subsequence, 
$$\int_{\Omega} \epsilon(\eu_{h}) \chi_{A_{h}}: \phi \mbox{ d}x 
\rightarrow \int_{\Omega} \epsilon : \phi \mbox{ d}\mu$$ for all 
$\phi \in C^{0}_{c} (\Omega, R^{n \times n}_{sym})$, where $(\mu, \epsilon)$ 
is a measure-function and $\mu \ll \leb$, $\frac{\mu}{\leb} \in [0,1]$ 
$\leb$ a.e. in $\Omega$. Also, 
$$\int_{\Omega} w(\epsilon) \mbox{ d}\mu \leq \liminf_{h \rightarrow \infty} 
\int_{\Omega} w(\epsilon(\eu_{h}) \mbox{ d}x$$ Suppose that $w$ is the 
well-known $C^{\infty}$ energy  potential of a linear hyperelastic material 
with the Hooke tensor $\eC$. The weak convergence of the sequence 
of measure-functions implies that for all $\ev \in C^{2}_{c}(\Omega \setminus 
K,R^{n})$ $$\int_{\Omega} \eC \epsilon(\ev) : \epsilon(\eu_{h}) \chi_{A_{h}} 
\mbox{ d}x \rightarrow \int_{\Omega} \eC \epsilon(\ev): \epsilon \mbox{ d}\mu$$ Recall that any 
$\eu_{h}$ is balanced and $\eC$ is a symmetric 4-order tensor;
 therefore for all $\ev \in C^{2}_{c}(\Omega \setminus K,R^{n})$ 
$$\int_{\Omega}\eC \epsilon : \epsilon(\ev) \mbox{ d}\mu = 0$$ which means that 
\begin{equation}
\frac{\mu}{\leb} \eC \epsilon \in Y_{K}(\Omega)
\label{bala}
\end{equation}  
For any $h \in N$ and for any $\phi \in C_{c}^{1}(\Omega,R^{n \times n}_{sym})$  the following 
equality holds:
$$\int_{\Omega} \eu_{h} div \phi \chi_{A_{h}} \mbox{ d}x + 
\int_{\Omega} \phi : \epsilon(\eu_{h}) \chi_{A_{h}} \mbox{ d}x + 
\int_{K \cap A_{h} } [\eu_{h}] \odot \en : \phi \mbox{ d}\hen = 0$$
Up to a subsequence $([\eu_{h}]\odot \en \chi_{A_{h}})_{h}$ converges weak * 
in $L^{\infty}$ to $a$ and the limit of the previous inequality is: for all 
$\phi \in C^{1}_{c}(\Omega,R^{n \times n})$
\begin{equation}
\int_{\Omega} \eu div \phi \mbox{ d}x + 
\int_{\Omega} \epsilon : \phi \mbox{ d}\mu + 
\int_{K} a : \phi \mbox{ d}\hen = 0
\label{inte}
\end{equation}
That justifies the notation $\epsilon = \epsilon_{\mu}(\eu)$. (\ref{bala}) shows nothing but the 
fact that $\eu$ is balanced.

This fact suggests to enlarge the domain of admissible triples in order to contain triples like 
$(\eu,\mu,K)$ which satisfy (\ref{inte}), $\mu \ll \leb$, $\frac{\mu}{\leb} \in 
[0,1]$. The enlarged energy functional is:
 $$ I(\eu,\mu,K) = \int_{\Omega} w(\epsilon_{\mu}(\eu)) - \gamma  \mbox{ d} \mu + 
G \hen(K)$$

We shall enlarge even more the set of admissible triples by allowing $K$ to be 
an integral varifold. With the same kind of argumentation as before we can prove that there exist 
a minimum of the functional
 $$I(\eu,\mu,\eV) = \int_{\Omega} w(\epsilon_{\mu}(\eu)) - 
\gamma \mbox{ d}\mu + G \| \eV \|(\Omega)$$ in the set of all admissible triples such that the first 
variation $\mid \delta \eV \mid$ is uniformly bounded. We do not repeat the proof 
because the use of the compactness theorem for integral varifolds is the only 
new element that appears. 

Denote by $X(\eu_{0},N)$ the closure of $T(\eu_{0},N)$ with respect to
the convergence:
$$(\eu_{h},\mu_{h},\eV_{h}) \rightarrow (\eu,\mu,\eV) \mbox{ iff}$$ 
$$ (\eu_{h},\mu_{h}) \rightarrow (\eu,\mu) \mbox{ as measure-functions}$$
$$\eV_{h} \rightarrow \eV \mbox{ as integral varifolds}$$
We have proved that for any sequence $(\eu_{h},\mu_{h},\eV_{h}) \in X(\eu_{0},N)$ such that 
$$\| \eu_{h} \|_{\infty} + I(\eu_{h},\mu_{h},\eV_{h})$$ is uniformly bounded, there is a 
subsequence converging to an element $(\eu,\mu,\eV)$ of a $X(\eu_{0},M)$ such that 
$$I(\eu,\mu,\eV) \leq \liminf_{h \rightarrow \infty} I(\eu_{h},\mu_{h},\eV_{h})$$
Also, for any $(\eu,A,K) \in T(\eu_{0},N)$  exists $M \in R$ and $(\ev,\mu,K) \in X(\eu_{0},N)$ with 
$$\| \ev \|_{\infty} + I(\ev,\mu,\eV) < + \infty$$ and $\ev$ balanced, i.e. 
for all $\ew \in C^{2}_{c}(\Omega \setminus K,R^{n})$ 
$$\int_{\Omega}\eC \epsilon_{\mu}(\ev) : \epsilon(\ew) \mbox{ d}\mu = 0 \mbox{ ,}$$  
such that 
\begin{equation}
\|w(\epsilon_{\mu}(\ev)) \frac{\mu}{\leb} \|_{\infty} \leq \gamma
\label{bound}
\end{equation} 
and $I(\eu,A,K) \geq I(\ev,\mu,K)$. In order to prove the existence of a minimizer of $I$ over 
$$X(\eu_{0}) = \bigcup_{N \in R} X(\eu_{0},N)$$ it is therefore sufficient to consider minimizing
 sequences formed by balanced terms. We use 
(\ref{firstb}) to show that  for $$N > \frac{\gamma}{G}$$  any minimizer 
 of $I$ over  $X(\eu_{0},M)$ belongs to $X(\eu_{0},N)$.  

Consider a sequence of balanced $(\eu_{h},A_{h},K_{h}) \in T(\eu_{0})$ converging in $X(\eu_{0})$ to 
$(\eu,\mu,\eV)$. Suppose more that (\ref{bound}) is true for all $(\eu_{h},A_{h},K_{h})$ and 
$I(\eu_{h},A_{h}K_{h})$ is uniformly bounded. Since 
$\delta_{\eta}I(\eu,\mu,\eV)$ can not be defined for any element of $X(\eu_{0})$, we adopt the 
following definition of the first variation of $I$: take a $\eta \in C^{2}_{0}(\Omega,R^{n})$,
$\mid \eta \mid \leq 1$, and 
denote by $\phi_{\epsilon}$ the one-group parameter generated by $\eta$. For any sequence 
$(\eu_{h},A_{h},K_{h})$ in 
$T(\eu_{0})$ converging to $(\eu,\mu,\eV)$ it is not hard to see that 
$(\eu_{h}.\phi_{\epsilon}^{-1} ,\phi_{\epsilon}(A_{h}),\phi_{\epsilon}(K_{h}))$ converges up to 
a subsequence to an element of $X(\eu_{0})$ denoted by $(\eu.\phi_{\epsilon}^{-1},\phi_{\epsilon}
(\mu),\phi_{\epsilon}(\eV))$ (attention, this element depends on the sequence 
$(\eu_{h},A_{h},K_{h})$). By definition
$$\delta_{\eta}I(\eu,\mu,\eV) = \inf \left\{ \liminf_{\epsilon \rightarrow 0} 
\frac{ I(\eu.\phi_{\epsilon}^{-1},\phi_{\epsilon}
(\mu),\phi_{\epsilon}(\eV)) - I(\eu,\mu,\eV)}{\mid \epsilon \mid} \mbox{ : } 
(\eu_{h},A_{h},K_{h}) \rightarrow (\eu,\mu,\eV) \right\}$$
(\ref{firstb}) and (\ref{bound}) show that $\delta_{\eta}I(\eu_{h},A_{h},K_{h})$ is uniformly 
bounded. For any $(\eu,A,K) \in T(\eu_{0})$ the following function is lower semicontinuous: 
$$\delta_{\eta}I(\eu,A,K) (\epsilon) = \left\{ 
\begin{array}{ll}
\frac{I(\eu.\phi_{\epsilon}^{-1},\phi_{\epsilon}(A),\phi_{\epsilon}(K)) - I(\eu,A,K)}{\mid 
\epsilon \mid} & \mbox{ if } \epsilon \not = 0 \\
\delta_{\eta}I(\eu,A,K) & \mbox{ if } \epsilon = 0
\end{array} \right. $$
Again the compactness theorem for measure-functions applied to the sequence $$(
\hen_{|_{K_{h}}}, 
[w(\epsilon(\eu_{h})] \en_{h})$$ associated to $(\eu_{h},A_{h},K_{h})$ affirms that 
(up to a subsequence) $( \hen_{|_{K_{h}}}, 
[w(\epsilon(\eu_{h})] \en_{h}) \rightarrow ( \| \eV \|, \lambda)$ with 
$\| \lambda \|_{\infty} \leq \gamma$. Because of the uniform upper 
bound of $\delta_{\eta}I(\eu_{h},A_{h},K_{h})$, the mentioned lower semicontinuity implies that 
$$\delta_{\eta}I(\eu,\mu,\eV) \leq   \int_{\Omega} \eta \cdot \lambda \mbox{ d}\| \eV \| 
+ G \delta \eV (\eta)$$ 
It is obvious that for any minimizer of $I$ over a $X(\eu_{0},M)$ 
$$\delta_{\eta} I(\eu,\mu,\eV) \geq 0$$
The latter two inequalities end the proof.

\vspace{.5cm}

{\bf Theorem 2:} {\it The functional $I$ has a minimizer over the set $X(\eu_{0})$. Moreover, 
any minimizer of $I$ over this set belongs to $X(\eu_{0}, 2 \gamma / G )$}

\vspace{.5cm}

The crack predicted by the second model (i.e. $\eV$) has no edges. Indeed, suppose that $\eV$ is 
a smooth manifold. Then the "cut" procedure eliminates the edge if it  belongs to the sound region. 
It follows that the edge of the crack lies in a damaged region. In other words, around the edge it 
is at least a region containing a fine mixture of damaged and undamaged material. 

\vspace{1.cm}

I want to thank to Prof. Luigi Ambrosio for his suggestions, especially concerning the $\eSBD$ 
approach to the first model.

\vspace{2.cm}

REFERENCES

\vspace{.5cm}

\begin{itemize}
\item
[Amb] L. Ambrosio, Variational problems in SBV and image segmentation, Acta Appl.  
Mathematic\ae  17, 1-40, 1989

\item
[AFB] L. Ambrosio, G. Buttazzo, I. Fonseca, Lower semicontinuity problems in Sobolev spaces 
with respect to a measure, J. Math. Pures Appl. 75, 1996

\item
[AFP] L. Ambrosio, N. Fusco, D. Pallara, Partial regularity of free discontinuity sets II, to 
appear

\item
[AP] L. Ambrosio, D. Pallara, Partial regularity of free discontinuity sets I, to appear

\item
[BCDM] G. Bellettini, A. Coscia, G. Dal Maso, Compactness and lower semicontinuity properties 
in $\eSBD(\Omega)$, preprint S.I.S.S.A. 86/96/M, 1996

\item
[Bra] K. A. Brakke, The motion of a surface by its mean curvature, Math. Notes, Princeton Univ. 
Press, 1978

\item
[DGCL] E. De Giorgi, G. Carriero, A.  Leaci, Existence theory for a minimum problem with free 
discontinuity set, Arch. Rational Mech. Anal., vol. 108, 1989

\item
[DP] G. Del Piero, Recent developments in the mechanics of materials which do not support tension, 
in Free Boundary Problems: Theory and 
Applications, vol I, Eds. Hoffmann K. H., Sprekels J., Pitman res. notes in math. series, 
Longman Scientific \& Technical, 1990

\item
[DS] G. David, S. Semmes, On the singular sets of minimizers of the Mumford-Shah functional, 
J. Math. Pures Appl. 75, 1996

\item
[FMa] G. Francfort, J.-J. Marigo, Stable damage evolution in a brittle continuous medium, 
Eur. J. Mech., A/Solids, 12, no. 2, 1993

\item
[FMu] G. Francfort, F. Murat, Homogenization and optimal bounds in linear elasticity, Arch. 
Rational Mech. Anal., 94, 1986

\item
[G1] M. E. Gurtin, On a Theory of Phase Transitions with Interfacial Energy, Arch. 
Rational Mech. Anal., 87, 1985

\item
[G2] M. E. Gurtin, On Phase Transitions with Bulk, Interfacial and Boundary Energy, Arch. 
Rational Mech. Anal., 96, 1986

\item
[GWZ] M. E. Gurtin, W. O. Williams \& W. P. Ziemer, Geometric Measure Theory and the Axioms of 
Continuum  Thermodynamics, Arch. Rational Mech. Anal., 92, 1986

\item
[Hut] J. E. Hutchinson, Second fundamental form for varifolds and the existence of surfaces 
minimizing curvature, Indiana Univ. Math. J., vol. 35, no. 1, 1986

\item
[M] J. J. Moreau, Fonctionnelles convexes, S\'eminaire sur les Equations aux Derivees Partielles, 
College de France 1966-1967

\item
[Mo] L. Modica, The Gradient Theory of Phase Transitions and the Minimal Interface Criterion, 
Arch. Rational Mech. Anal., 98, 1987 

\item
[MS] D. Mumford, J. Shah, Optimal approximation by piecewise smooth functions and associated 
variational problems, Comm. on Pure and Appl. Math., vol. XLII, 
no. 5, 1989

\end{itemize}

\end{document}